\documentclass{amsart}

\usepackage{mathtools}
\usepackage{amsmath}
\usepackage{amssymb}
\usepackage{graphicx}
\usepackage{amsthm}
\usepackage[square,numbers]{natbib}
\usepackage{tikz}

\DeclareMathOperator*{\DDottedprod}{%
\begin{tikzpicture}[baseline=(T.base)]
\node[inner sep=0pt, outer sep=0pt] (T) at (0,0) {$\displaystyle\prod$};
\node[inner sep=0pt, outer sep=0pt] (T) at (0,0) {$\displaystyle\coprod$};
\draw[fill] (0,2.5pt) circle (1.5pt);
\draw[fill] (0,-2.5pt) circle (1.5pt);
\end{tikzpicture}}

\newtheorem{theorem}{Theorem}
\newtheorem{lemma}{Lemma}
\newtheorem{corollary}{Corollary}[theorem]

\theoremstyle{definition}

\author{Efe Gürel}
\address{TÜBİTAK Natural Sciences High School, Kocaeli, 41400, Turkey}
\email{efegurel54@gmail.com}

\title[Zeta Regularized Trigonometric Products]{Zeta Regularized Trigonometric Products Over Zeros Of The Riemann Zeta Function}


\begin{document}
\begin{abstract}
    We prove a novel zeta regularized product formula concerning regularization of trigonometric products over non-trivial zeros of the Riemann zeta function. Furthermore, we calculate the discrepancies of such regularized products. In special cases, our formula reduces to the Kimoto-Wakayama formula. A conjectural relationship between such products and a weak Riemann hypothesis is speculated.
\end{abstract}
\maketitle

\section{Introduction}
Let $\Lambda=\left\{ \lambda_n \right\}_{n\in I}$ be a sequence of non-zero complex numbers. Associated zeta function of $\Lambda$, denoted $\zeta_\Lambda$, is defined as
\begin{align*}
    \zeta_\Lambda(s)=\sum_{n\in I}\frac{1}{\lambda_n^s}
\end{align*}
which we assume to be absolutely convergent for sufficiently large $\mathfrak{Re}(s)$. Here and throughout, the principal branch of the logarithm is chosen. Suppose that $\zeta_\Lambda(s)$ is analytically continued to $\mathfrak{Re}(s)>0$ and there exists a finite sequence of functions $\left\{ f_k(s;\Lambda) \right\}_{k=1}^{n}$, each meromorphic at $s=0$, such that the function 
\begin{align} \label{RegularizableDef}
    P(s;\Lambda)=\zeta_\Lambda(s)-\sum_{k=1}^{n}f_k(s;\Lambda)(\log s)^k,
\end{align}
is analytically continued to a region containing $s=0$ as a single valued meromorphic function. Then we say that function $\zeta_\Lambda$ and sequence $\Lambda$ are regularizable. The function $P(s;\Lambda)$ is called the meromorphic part of $\zeta_\Lambda(s)$ at $s=0$. For a regularizable zeta function $\zeta_\Lambda$, the linear term of $\zeta_\Lambda$ at $s=0$ is defined as
\begin{align*}
    \underset{s=0}{\mathcal{LT}}\ \zeta_\Lambda(s)=\underset{s=0}{\text{Res}}\frac{P(s;\Lambda)}{s^2}.
\end{align*}
For a regularizable sequence $\Lambda$, the ddotted product of $\Lambda$ \cite{WakayamaRemarks} is defined by  

\begin{align*}
\DDottedprod_{n\in I}\lambda_n=\exp\left(-\underset{s=0}{\mathcal{LT}}\ \zeta_\Lambda(s)\right).   
\end{align*}
Although regularized products share a vast majority of properties related to classical products, it is not always true that $\DDottedprod \lambda_n\mu_n=\DDottedprod \lambda_n \ \DDottedprod\mu_n$. Let $\Lambda_k=\left\{ \lambda_{n,k} \right\}_{n\in I}$ for $k=1,\ldots,m$ and their pointwise product $\Lambda=\left\{ \lambda_{n,1}\lambda_{n,2}\cdots \lambda_{n,m} \right\}_{n\in I}$ be regularizable sequences. The discrepancy of the sequences $\Lambda_k$ is defined by
\begin{align*}
    F=F(\Lambda_1,\Lambda_2,\ldots,\Lambda_m)=\log\frac{\DDottedprod_{n\in I}\prod_{k=1}^{m}\lambda_{n,k}}{\prod_{k=1}^{m}\DDottedprod_{n\in I}\lambda_{n,k}}=\sum_{k=1}^{m}\underset{s=0}{\mathcal{LT}}\ \zeta_{\Lambda_k}(s)-\underset{s=0}{\mathcal{LT}}\ \zeta_\Lambda(s).
\end{align*}
Existence of this discrepancy is called the multiplicative anomaly. The problem of determining the discrepancy between regularized products is the main goal of theory of zeta regularized products. We refer the reader to \cite{Discrepancies,FiniteProducts,Resolution,BasicAnalysis,Mizuno,Quine,Voros} and references therein for discussion of zeta regularized products and multiplicative anomaly.
\\

The Riemann zeta function is given by
\begin{align*}
    \zeta(s)=\sum_{n=1}^{\infty}\frac{1}{n^s} \qquad (\mathfrak{Re}(s)>1)
\end{align*}
and can be analytically continued to the whole complex plane except $s=1$ where $\zeta(s)$ has a simple pole with residue $1$ \cite{TheoryofZeta}. The non-trivial zeros of Riemann zeta function are denoted by $\rho$ and $\rho=1/2+i\tau$. Cram{\'e}r's $V$ and $\phi$ functions \cite{Cramer} are defined as
\begin{align*}
    \begin{split}
        V(s)=\sum_{\mathfrak{Im}(\rho)>0}e^{s\rho} \qquad(\mathfrak{Im}(s)>0), \\
        \phi(s)=\sum_{\mathfrak{Re}(\tau)>0}e^{-s\tau} \qquad(\mathfrak{Re}(s)>0),
    \end{split}
\end{align*}
where summations are taken over all non-trivial zeros of $\zeta(s)$ under given conditions. Let $\gamma$ denote the Euler's constant. Cram{\'e}r also proved that the function
\begin{align} \label{CramerResult}
    V(s)-\frac{1}{2\pi i}\left(\frac{\log s}{1-e^{-s}}+\frac{\gamma+\log2\pi-\pi i/2}{s}\right)
\end{align}
can be analytically continued to whole complex plane and is holomorphic near $s=0$.
Two examples of trigonometric ddotted products over $\rho$ are given in \cite{WakayamaRemarks} as the following theorems.
\begin{theorem} \label{firstthm}
    For $\alpha>0$ and $x\in \mathbb{C}$ such that $2\pi/\alpha>\mathfrak{Re}(x)\ge 0$ and $\mathfrak{Im}(x)\le 0$, the function defined as the ddotted product
    \begin{align*}
        S_\alpha(x)=\DDottedprod_{\mathfrak{Im}(\rho)>0}\sin\alpha(\rho-x)
    \end{align*}
    exists. Furthermore,
    \begin{align*}
        \begin{split}
            S_\alpha(x)&=e^{-F_\alpha(x)}\prod_{\mathfrak{Im}(\rho)>0}e^{i\alpha(\rho-x)+\log-2i}\sin\alpha(\rho-x) \\ &=e^{-F_\alpha(x)}\left(e^{-2i \alpha x};e^{-2i\alpha }\right)_\zeta
        \end{split}
    \end{align*}
    where $(x;q)_\zeta$ is the zeta Pochhammer symbol defined as
    \begin{align*}
        (x;q)_\zeta=\prod_{\mathfrak{Im}(\rho)>0}\left(1-xq^{-\rho}\right)
    \end{align*}
    and $F_\alpha(x)$ is a quadratic polynomial in $x$ such that $F_\alpha(x)=A_\alpha x^2+B_\alpha x+C_\alpha$ and constants $A_\alpha$ and $B_\alpha$ are given by
    \begin{align*}
        A_\alpha=\frac{\alpha(\gamma+\log 2\pi\alpha)}{4\pi}, \qquad B_\alpha=\frac{(\gamma+\log 2\pi\alpha)}{2\pi i}\left(\frac{i\alpha}{2}+\log-2i\right)+\frac{7i\alpha}{8}.
    \end{align*}
\end{theorem}
\begin{theorem}
    For every $x\in \mathbb{C}$ and some constant $c$, the following equality holds.
    \begin{align*}
        \prod_{\mathfrak{Im}(\rho)>0}e^{i(x-\rho)}=\exp\left(-\frac{\gamma+\log2\pi}{4\pi}\left(x-\frac{1}{2}\right)^2-\frac{7ix}{8}+c\right)
    \end{align*}
\end{theorem}
In this paper, we consider functions defined as the trigonometric ddotted products of forms
\begin{align*}
    \DDottedprod_{\mathfrak{Im}(\rho)>0}\prod_{k=1}^{n}\sin (\alpha_k\rho-z_k)
\end{align*}
and
\begin{align*}
     \qquad \DDottedprod_{\mathfrak{Im}(\rho)>0}\prod_{k=1}^{n}\left(e^{-i(\alpha_k\rho-z_k)}-\omega_k\right)
\end{align*}
for suitable conditions on $\alpha_k,z_k$ and $\omega_k$. Furthermore, we calculate the discrepancies associated to such sequences.
\section{Main Results}

Throughout this section, let $\left\{ \alpha_k \right\}_{k=1}^n$ and $\left\{ z_k \right\}_{k=1}^n$ be sequences of complex numbers such that for every $k=1,\ldots,n$, inequalities $2\pi>\mathfrak{Re}(z_k)\ge0\ge\mathfrak{Im}(z_k)$ and $\alpha_k>0$ hold. Furhthermore, let us denote $z=z_1+z_2+...+z_n$ and $\alpha=\alpha_1+\alpha_2+...+\alpha_n$. We define the asssociated zeta function on the sequence $\prod_{k=1}^{n}\sin (\alpha_k\rho-z_k)$ as
\begin{align} \label{LSeriesDef}
    L(s;\underline{z},\underline{\alpha})=\sum_{\mathfrak{Im}(\rho)>0}\prod_{k=1}^{n}(\sin \left(\alpha_k\rho-z_k)\right)^{-s}
\end{align}
where $\underline{z}=(z_1,z_2,...,z_n)$ and $\underline{\alpha}=(\alpha_1,\alpha_2,...,\alpha_n)$. We shall explain the usage of the above zeta function instead of the associated one
\begin{align*}
    \ell(s;\underline{z},\underline{\alpha})=\sum_{\mathfrak{Im}(\rho)>0}\left( \prod_{k=1}^{n}\sin(\alpha_k\rho-z_k) \right)^{-s}.
\end{align*} 
It has been noted in \cite{Mizuno} that taking a different branch of logarithm for finitely many terms in the zeta function does not affect the value of the regularized product. Therefore, two zeta functions will result in the same regularization process if the equation 
\begin{align*}
    \prod_{k=1}^{n}\sin(\alpha_k\rho-z_k)^{-s}=\left( \prod_{k=1}^{n}\sin(\alpha_k\rho-z_k) \right)^{-s}
\end{align*}
holds for all but finitely many $\rho$. We know that this is true if we have
\begin{align*}
    -\pi < \sum_{k=1}^{n}\arg\sin(\alpha_k\rho-z_k) \le \pi.
\end{align*}
As $\mathfrak{Im}(\rho)\to \infty$, we have the asymptotic relation
\begin{align*}
    \arg\sin(\alpha_k\rho-z_k)&\sim\arg(\sin(\alpha_k\mathfrak{Re}(\rho)-\mathfrak{Re}(z_k))+i\cos (\alpha_k\mathfrak{Re}(\rho)-\mathfrak{Re}(z_k))) \\
    &=\arg ie^{-i(\alpha_k\mathfrak{Re}(\rho)-\mathfrak{Re}(z_k))}\\
    &=\frac{\pi}{2}-(\alpha_k\mathfrak{Re}(\rho)-\mathfrak{Re}(z_k)) \bmod 2\pi.
\end{align*}
Thus if one assumes the condition
\begin{align}\label{ArgCond1}
        -\pi < \sum_{k=1}^{n}\arg ie^{-i(\alpha_k\mathfrak{Re}(\rho)-\mathfrak{Re}(z_k))}< \pi.
\end{align}
holds true for all but finitely many $\rho$, the functions $ L(s;\underline{z},\underline{\alpha})$ and $ \ell(s;\underline{z},\underline{\alpha})$ will yield the same regularized products. We remark that this inequality is not heavily dependent on $\rho$ as $0<\mathfrak{Re}(\rho)<1$ and is conjectured to be always $1/2$. Our first main result concerns the regularizability of $ L(s;\underline{z},\underline{\alpha})$.
\begin{theorem}
    Assume that the inequality \eqref{ArgCond1} holds for all but finitely many $\rho$. Then, the function $L(s;\underline{z},\underline{\alpha})$ is regularizable.
\end{theorem}
\begin{proof}
    Replacing $\sin z=\frac{e^{iz}-e^{-iz}}{2i}$ in the equation \eqref{LSeriesDef}, we obtain
    \begin{align} \label{LSeriesToSine}
        L(s;\underline{z},\underline{\alpha})=\sum_{\mathfrak{Im}(\rho)>0}\prod_{k=1}^{n}(-2i)^se^{is(\alpha_k\rho-z_k)}\left(1-e^{2i(\alpha_k\rho-z_k)}\right)^{-s}.
    \end{align}
    Rearranging equation \eqref{LSeriesToSine} yields
    \begin{align} \label{LSeriesBeforeBinom}
        L(s;\underline{z},\underline{a})=(-2i)^{ns}e^{-isz}\sum_{\mathfrak{Im}(\rho)>0}e^{i\alpha s\rho}\prod_{k=1}^{n}\left(1-e^{2i(\alpha_k\rho-z_k)}\right)^{-s}.
    \end{align}
    Since $\left| e^{2i(\alpha_k\rho-z_k)} \right|=e^{2\mathfrak{Im}(z)-2\alpha_k\mathfrak{Re}(\tau)}<1$, we can apply the binomial theorem to get 
    \begin{align} \label{BinomialSeries}
        \left(1-e^{2i(\alpha_k\rho-z_k)}\right)^{-s}=\sum_{m=0}^{\infty}(-1)^m\binom{-s}{m}e^{2im(\alpha_k\rho-z_k)}.
    \end{align}
    Therefore by substituting \eqref{BinomialSeries} in equation \eqref{LSeriesToSine}, for certain constants $P_m(\underline{z},\underline{\alpha})$ and $Q_m(\underline{z},\underline{\alpha})$, we have that,
    \begin{align} \label{LSeriesRegu1}
         \begin{split}
            L(s;\underline{z},\underline{\alpha})&=(-2i)^{ns}e^{-isz}\sum_{\mathfrak{Im}(\rho)>0}e^{i\alpha s\rho}\prod_{k=1}^{n}\sum_{m=0}^{\infty}(-1)^m\binom{-s}{m}e^{2im(\alpha_k\rho-z_k)}\\
            &=(-2i)^{ns}e^{-isz}\sum_{\mathfrak{Im}(\rho)>0}e^{i\alpha s\rho}\sum_{m=0}^{\infty}s^m P_m(\underline{z},\underline{\alpha})e^{Q_m(\underline{z},\underline{\alpha})\rho}.
         \end{split}    
    \end{align}
   Interchanging order of summations in the equation \eqref{LSeriesRegu1} results 
    \begin{align} \label{LSeriesRegu2}
        \begin{split}
            L(s;\underline{z},\underline{\alpha})
            &=(-2i)^{ns}e^{-isz}\sum_{m=0}^{\infty}s^mP_m(\underline{z},\underline{\alpha})\sum_{\mathfrak{Im}(\rho)>0}e^{(Q_m(\underline{z},\underline{\alpha})+i\alpha s)\rho} \\ 
            &=(-2i)^{ns}e^{-isz}\sum_{m=0}^{\infty}s^mP_m(\underline{z},\underline{\alpha})V(Q_m(\underline{z},\underline{\alpha})+i\alpha s).
        \end{split}
    \end{align}
    Hence the function $L(s;\underline{z},\underline{\alpha})$ is regularizable by the equations \eqref{CramerResult} and \eqref{LSeriesRegu2} combined with \eqref{RegularizableDef}. Thus completing the proof.
\end{proof}
To find the linear term of $L(s;\underline{z},\underline{\alpha})$, we need the following lemma given in \cite{WakayamaRemarks}.

\begin{lemma} \label{MeroPhiLemma}
    For $\alpha>0$, the meromorphic part $\varphi_\alpha(s)$ of the function $\phi(\alpha s)$ is given by
    \begin{align*}
        \varphi_\alpha(s)=c_{-1}s^{-1}+c_0+c_1s+O\left(s^2\right)
    \end{align*}
    where coefficients $c_{-1}$ and $c_0$ satisfy
    \begin{align*}
        c_{-1}=-\frac{\gamma+\log2\pi \alpha}{2\pi \alpha}, \qquad c_0=\frac{7}{8}.
    \end{align*}
\end{lemma}

\begin{theorem}\label{MainThm1}
    Assume that the inequality \eqref{ArgCond1} holds for all but finitely many $\rho$. Then, the function defined as the ddotted product 
    \begin{align*}
             S(\underline{z};\underline{\alpha})=\DDottedprod_{\mathfrak{Im}(\rho)>0}\prod_{k=1}^{n}\sin (\alpha_k\rho-z_k)
    \end{align*}
    exists. Furthermore,
        \begin{align*}
            \begin{split}
                S(\underline{z};\underline{\alpha})&=e^{-F(\underline{z};\underline{\alpha})}\prod_{k=1}^{n}\prod_{\mathfrak{Im}(\rho)>0}\sin(\alpha_k\rho-z_k)e^{i(\alpha_k\rho-z_k)+\log-2i} \\
                &=e^{-F(\underline{z};\underline{\alpha})}\prod_{k=1}^{n}\left(e^{-2iz_k};e^{-2i\alpha_k}\right)_\zeta
            \end{split}
        \end{align*}
        where $F(\underline{z};\underline{\alpha})$ is a quadratic polynomial given by
        \begin{align*}
            F(\underline{z};\underline{\alpha})=&-\frac{(\gamma+\log2\pi \alpha)}{4\pi \alpha}\left(n\log(-2i)-iz+\frac{i\alpha}{2}\right)^2\\
            &+\frac{7(2n\log(-2i)-2iz+i\alpha)}{16}+c_1.
        \end{align*}
\end{theorem}
\begin{proof}
    Starting with equation \eqref{LSeriesBeforeBinom} and \eqref{BinomialSeries}, we have that 
    \begin{align*}
        L(s;\underline{z},\underline{\alpha})=(-2i)^{ns}e^{-isz}\sum_{\mathfrak{Im}(\rho)>0}e^{i\alpha s\rho}\prod_{k=1}^{n}\sum_{m=0}^{\infty}(-1)^m\binom{-s}{m}e^{2im(\alpha_k\rho-z_k)}.
    \end{align*}
    Isolating the term $s$ in every binomial coefficient, we have $(-1)^m\binom{-s}{m}=\frac{s}{m}+O\left(s^2\right)$ for $m\ge1$, which means
    \begin{align*}
    \sum_{m=0}^{\infty}(-1)^m\binom{-s}{m}e^{2im(\alpha_k\rho-z_k)}=1+s\sum_{m=1}^{\infty}\frac{e^{2im(\alpha_k\rho-z_k)}}{m}+O\left(s^2\right)
    \end{align*}
    and therefore we obtain
    \begin{align} \label{LSeriesAfterBinomMultiplication}
        \begin{split}
            L(s;\underline{z},\underline{\alpha})&=(-2i)^{ns}e^{-isz}\sum_{\mathfrak{Im}(\rho)>0}e^{i\alpha s\rho}\prod_{k=1}^{n}\left(1+s\sum_{m=1}^{\infty}\frac{e^{2im(\alpha_k\rho-z_k)}}{m}+O\left
            (s^2\right)\right)\\
            &=(-2i)^{ns}e^{-isz}\sum_{\mathfrak{Im}(\rho)>0}e^{i\alpha s\rho}\left(1+s\sum_{k=1}^{n}\sum_{m=1}^{\infty}\frac{e^{2im(\alpha_k\rho-z_k)}}{m}+O\left(s^2\right)\right).
        \end{split}
    \end{align}
    Distributing the summation in equation \eqref{LSeriesAfterBinomMultiplication} yields
    \begin{align} \label{LSeriesToVFunction}
        \begin{split}
            L(s;\underline{z},\underline{\alpha})&=(-2i)^{ns}e^{-isz}\left(V(i\alpha s)+s\sum_{k=1}^{n}\sum_{m=1}^{\infty}\sum_{\mathfrak{Im}(\rho)>0}\frac{e^{(2im\alpha_k+i\alpha s)\rho-2imz_k}}{m}+O\left(s^2\right)\right) \\
            &=e^{s(-iz+n\log-2i)}\left(V(i\alpha s)+s\sum_{k=1}^{n}\sum_{m=1}^{\infty}\frac{V(2im\alpha_k+i\alpha s)e^{-2imz_k}}{m}+O\left(s^2\right)\right).
        \end{split}
    \end{align}
    Since the function $V(2im\alpha_k+i\alpha s)$ is holomorphic at $s=0$, we have $V(2im\alpha_k+i\alpha s)=V(2im\alpha_k)+O(s)$ and
    \begin{align} \label{ExpandingV}
        \sum_{m=1}^{\infty}\frac{V(2im\alpha_k+i\alpha s)e^{-2imz_k}}{m}=\sum_{m=1}^{\infty}\frac{V(2im\alpha_k)e^{-2imz_k}}{m}+O(s).
    \end{align}
    For brevity, defining the function 
    \begin{align} \label{AuxFunctionDef}
        f(k)=\sum_{m=1}^{\infty}\frac{V(2im\alpha_k)e^{-2imz_k}}{m}
    \end{align}
    and substituting equation \eqref{ExpandingV} in \eqref{LSeriesToVFunction}, we get
    \begin{align} \label{LClosedFormFinally}
        L(s;\underline{z},\underline{\alpha})=e^{s(-iz+n\log-2i)}\left(V(i\alpha s)+s\sum_{k=1}^{n}f(k)+O\left(s^2\right)\right).
    \end{align}
    Taking linear terms in both sides of equation \eqref{LClosedFormFinally} results
    \begin{align*} 
        \begin{split}
            \underset{s=0}{\mathcal{LT}}\ L(s;\underline{z},\underline{\alpha})&=\underset{s=0}{\mathcal{LT}} \left[ e^{s(-iz+n\log-2i)}\left(V(i \alpha s)+s\sum_{k=1}^{n}f(k)+O\left(s^2\right)\right)\right] \\
            &=\underset{s=0}{\mathcal{LT}}\left(e^{s(-iz+n\log-2i+i\alpha/2)}\phi(as)\right)+\sum_{k=1}^{n}f(k)
        \end{split}
    \end{align*}
    where we have used the trivial fact that $V(iz)=e^{iz/2}\phi(z)$. Let $A=-iz+n\log-2i+i\alpha/2$, by Lemma \ref{MeroPhiLemma}, we obtain
    \begin{align*}
        \underset{s=0}{\mathcal{LT}}\left(e^{sA}\phi(\alpha s)\right)=\underset{s=0}{\text{Res}}\frac{e^{sA}\varphi_\alpha(s)}{s^2}=\frac{A^2c_{-1}}{2}+Ac_0+c_1=:F(\underline{z};\underline{\alpha}).
    \end{align*}
    By the definition of ddotted product, we have
    \begin{align} \label{LinearTermToDDotted}
    \begin{split}
        S(\underline{z};\underline{\alpha})=\exp\left(-\underset{s=0}{\mathcal{LT}}\ L(s;\underline{z};\underline{\alpha})\right)&=\exp\left(-F(\underline{z};\underline{\alpha})-\sum_{k=1}^{n}f(k)\right)\\
        &=e^{-F(\underline{z};\underline{\alpha})}\prod_{k=1}^{n}\exp\left(-f(k)\right)
    \end{split}
    \end{align}
    Taking exponents in both sides of equation \eqref{AuxFunctionDef}, we get
    \begin{align*}
        \exp(-f(k))&=\exp\left(-\sum_{m=1}^{\infty}\frac{V(2im\alpha_k)e^{-2imz_k}}{m}\right)\\
        &=\exp\left(-\sum_{m=1}^{\infty}\sum_{\mathfrak{Im}(\rho)>0}\frac{e^{2im(\alpha_k\rho-z_k)}}{m}\right).
    \end{align*}
    We now interchange order of summations as follows,
    \begin{align} \label{Expf(k)}
        \begin{split}
            \exp(-f(k))&=\exp\left(-\sum_{\mathfrak{Im}(\rho)>0}\sum_{m=1}^{\infty}\frac{e^{2im(\alpha_k\rho-z_k)}}{m}\right)\\
            &=\exp\left(\sum_{\mathfrak{Im}(\rho)>0}\log\left(1-e^{2i(\alpha_k\rho-z_k)}\right)\right) \\
            &=\prod_{\mathfrak{Im}(\rho)>0}1-e^{2i(\alpha_k\rho-z_k)}\\
            &=\prod_{\mathfrak{Im}(\rho)>0}-2ie^{i\left(\alpha_k\rho-z_k\right)}\sin(\alpha_k\rho-z_k).
        \end{split}
    \end{align}
    Finally, substituting equation \eqref{Expf(k)} in \eqref{LinearTermToDDotted}, we have
    \begin{align*}
    \begin{split}
        S(\underline{z};\underline{\alpha})&=e^{-F(\underline{z};\underline{\alpha})}\prod_{k=1}^{n}\prod_{\mathfrak{Im}(\rho)>0}\sin(\alpha_k\rho-z_k)e^{i(\alpha_k\rho-z_k)+\log-2i}\\
        &=e^{-F(\underline{z};\underline{\alpha})}\prod_{k=1}^{n}\left(e^{-2iz_k};e^{-2i\alpha_k}\right)_\zeta.
    \end{split}
    \end{align*}
    This completes the proof.
\end{proof}
The discrepancy of the sequences $\sin(\alpha_k\rho-z_k)$ can be calculated easily by the means of Theorem \ref{MainThm1}.
\begin{corollary}
    The discrepancy $F$ of the sequences $\sin(\alpha_k\rho-z_k)$, $k=1,\ldots,n$ is given by
    \begin{align*}
        F=\sum_{k=1}^{m}F(z_k;\alpha_k)-F(\underline{z};\underline{\alpha}).
    \end{align*}
\end{corollary}
From this result, the discrepancy between sequences of form $\prod_{j=1}^{m_k}\sin(\alpha_{j,k}\rho-z_{j,k})$ may also be calculated inductively. Now we consider regularized products of the form
\begin{align*}
         \qquad \DDottedprod_{\mathfrak{Im}(\rho)>0}\prod_{k=1}^{n}\left(e^{-i(\alpha_k\rho-z_k)}-\omega_k\right).
\end{align*}
Similarly, we define the associated zeta function of the sequence $\prod_{k=1}^{n}\left(e^{-i(\alpha_k\rho-z_k)}-\omega_k\right)$ by
\begin{align} \label{L2SeriesDef}
    \widetilde{L}(s;\underline{z},\underline{\alpha},\underline{\omega})=\sum_{\mathfrak{Im}(\rho)>0}\prod_{k=1}^{n}\left(e^{-i(\alpha_k\rho-z_k)}-\omega_k\right)^{-s}
\end{align}
where $1\ge \left| \omega_k \right|$ holds for every $k$ and $\underline{\omega}=(\omega_1,\omega_2,...,\omega_n)$. By a similar comparison with the zeta function 
\begin{align*}
     \widetilde{\ell}(s;\underline{z},\underline{\alpha},\underline{\omega})=\sum_{\mathfrak{Im}(\rho)>0}\left(\prod_{k=1}^{n}\left(e^{-i(\alpha_k\rho-z_k)}-\omega_k\right)\right)^{-s},
\end{align*}
we see that these two zeta functions will result in the same regularization process if one assumes the condition
\begin{align} \label{ArgCond2}
            -\pi < \sum_{k=1}^{n}\arg e^{-i(\alpha_k\mathfrak{Re}(\rho)-\mathfrak{Re}(z_k))}< \pi
\end{align}
holds true for all but finitely many $\rho$. Again, we first prove the regularizability of $\widetilde{L}(s;\underline{z},\underline{\alpha},\underline{\omega})$. 
\begin{theorem}
    Assume that the inequality \eqref{ArgCond2} holds for all but finitely many $\rho$. Then, the function $\widetilde{L}(s;\underline{z},\underline{\alpha},\underline{\omega})$ is regularizable.
\end{theorem}
\begin{proof}
    Starting with equation \eqref{L2SeriesDef} and rearranging, we obtain 
    \begin{align} \label{L2Series1}
    \begin{split}
        \widetilde{L}(s;\underline{z},\underline{\alpha},\underline{\omega})&=\sum_{\mathfrak{Im}(\rho)>0}\prod_{k=1}^{n}e^{is(\alpha_k\rho-z_k)}\left(1-\omega_ke^{i(\alpha_k\rho-z_k)}\right)^{-s}\\
        &=e^{-isz}\sum_{\mathfrak{Im}(\rho)>0}e^{i\alpha s\rho}\prod_{k=1}^{n}\left(1-\omega_ke^{i(\alpha_k\rho-z_k)}\right)^{-s}.
    \end{split}
    \end{align}
    Since $\left| \omega_ke^{i(\alpha_k\rho-z_k)} \right|=|\omega_k|e^{\mathfrak{Im}(z)-\alpha_k\mathfrak{Re}(\tau)}<1$, we can use the binomial theorem to get 
    \begin{align} \label{BinomialSeriesW}
        \left(1-\omega_ke^{i(\alpha_k\rho-z_k)}\right)^{-s}=\sum_{m=0}^{\infty}(-1)^m\binom{-s}{m}\omega_k^me^{im(\alpha_k\rho-z_k)}.
    \end{align}
    Therefore by substituting equation \eqref{BinomialSeriesW} in equation \eqref{L2Series1}, for certain constants $\widetilde{P}_m(\underline{z};\underline{\alpha},\underline{\omega})$  and $\widetilde{Q}_m(\underline{z};\underline{\alpha},\underline{\omega})$, we have that
    \begin{align} \label{L2SeriesRegu1}
        \begin{split}
            \widetilde{L}(s;\underline{z},\underline{\alpha},\underline{\omega})&=e^{-isz}\sum_{\mathfrak{Im}(\rho)>0}e^{i\alpha s\rho}\prod_{k=1}^{n}\sum_{m=0}^{\infty}(-1)^m\binom{-s}{m}\omega_k^me^{im(\alpha_k\rho-z_k)}\\
            &=e^{-isz}\sum_{\mathfrak{Im}(\rho)>0}e^{i\alpha s\rho}\sum_{m=0}^{\infty}s^m\widetilde{P}_m(\underline{z};\underline{\alpha},\underline{\omega})e^{\widetilde{Q}_m(\underline{z};\underline{\alpha},\underline{\omega})\rho}.
        \end{split}
    \end{align}
    Interchanging orders of summation in equation \eqref{L2SeriesRegu1} yields 
    \begin{align} \label{L2SeriesRegu2}
        \begin{split}
            \widetilde{L}(s;\underline{z},\underline{\alpha},\underline{\omega})&=e^{-isz}\sum_{m=0}^{\infty}s^m\widetilde{P}_m(\underline{z};\underline{\alpha},\underline{\omega})\sum_{\mathfrak{Im}(\rho)>0}e^{\left(\widetilde{Q}_m(\underline{z};\underline{\alpha},\underline{\omega})+i \alpha s\right)\rho}\\
            &=e^{-isz}\sum_{m=0}^{\infty}s^m\widetilde{P}_m(\underline{z};\underline{\alpha},\underline{\omega})V\left(\widetilde{Q}_m(\underline{z};\underline{\alpha},\underline{\omega})+i\alpha s\right).
        \end{split}
    \end{align}
    Hence the function $\widetilde{L}(s;\underline{z},\underline{\alpha},\underline{\omega})$ is regularizable by equations \eqref{CramerResult} and \eqref{L2SeriesRegu2} combined with \eqref{RegularizableDef}. Thus completing the proof.
\end{proof}
\begin{theorem}\label{Mainthm2}
    Assume that the inequality \eqref{ArgCond2} holds for all but finitely many $\rho$. Then, the function defined as the ddotted product
     \begin{align*}
         \widetilde{S}(\underline{z};\underline{\alpha},\underline{\omega})=\DDottedprod_{\mathfrak{Im}(\rho)>0}\prod_{k=1}^{n}\left(e^{-i(\alpha_k\rho-z_k)}-\omega_k\right)
     \end{align*}
     exists. Furthermore
     \begin{align*}
         \begin{split}
             \widetilde{S}(\underline{z};\underline{\alpha},\underline{\omega})&=e^{-\widetilde{F}(\underline{z};\underline{\alpha})}\prod_{k=1}^{n}\prod_{\mathfrak{Im}(\rho)>0}1-\omega_ke^{i(\alpha_k\rho-z_k)}\\
             &=e^{-\widetilde{F}(\underline{z};\underline{\alpha})}\prod_{k=1}^{n}\left(\omega_ke^{-iz_k};e^{-i\alpha_k}\right)_\zeta
         \end{split}
     \end{align*}
     where $\widetilde{F}(\underline{z};\underline{\alpha})$ is a polynomial given by
     \begin{align*}
         \widetilde{F}(\underline{z};\underline{\alpha})=-\frac{\gamma+\log2\pi \alpha}{4\pi \alpha}\left(\frac{i\alpha}{2}-iz\right)^2+\frac{7(i\alpha-2iz)}{16}+c_1.
     \end{align*}
\end{theorem}
\begin{proof}
    Starting with equation \eqref{BinomialSeriesW} and equation \eqref{L2Series1}, we have that
    \begin{align*}
        \widetilde{L}(s;\underline{z},\underline{\alpha},\underline{\omega})=e^{-isz}\sum_{\mathfrak{Im}(\rho)>0}e^{i\alpha s\rho}\prod_{k=1}^{n}\sum_{m=0}^{\infty}(-1)^m\binom{-s}{m}\omega_k^me^{im(\alpha_k\rho-z_k)}.
    \end{align*}
    Similarly, isolating term $s$ in every binomial coefficient, we get
    \begin{align*}
        \sum_{m=0}^{\infty}(-1)^m\binom{-s}{m}\omega_k^me^{im(\alpha_k\rho-z_k)}=1+s\sum_{m=1}^{\infty}\frac{\omega_k^me^{im(\alpha_k\rho-z_k)}}{m}+O\left(s^2\right)
    \end{align*}
    and therefore we obtain
    \begin{align} \label{L2SeriesAfterBinomMultiplication}
        \begin{split}
            \widetilde{L}(s;\underline{z},\underline{\alpha},\underline{\omega})&=e^{-isz}\sum_{\mathfrak{Im}(\rho)>0}e^{i\alpha s\rho}\prod_{k=1}^{n}\left(1+s\sum_{m=1}^{\infty}\frac{\omega_k^me^{im(\alpha_k\rho-z_k)}}{m}+O\left(s^2\right)\right) \\ 
            &=e^{-isz}\sum_{\mathfrak{Im}(\rho)>0}e^{i\alpha s\rho}\left(1+s\sum_{k=1}^{n}\sum_{m=1}^{\infty}\frac{\omega_k^me^{im(\alpha_k\rho-z_k)}}{m}+O\left(s^2\right)\right).
        \end{split}
    \end{align}
    Distributing the summation in equation \eqref{L2SeriesAfterBinomMultiplication} yields
    \begin{align} \label{L2SeriesToVfunction}
        \begin{split}
        \widetilde{L}(s;\underline{z},\underline{\alpha},\underline{\omega})&=e^{-isz}\left(V(i\alpha s)+s\sum_{k=1}^{n}\sum_{m=1}^{\infty}\frac{\omega_k^me^{-imz_k}}{m}\sum_{\mathfrak{Im}(\rho)>0}e^{(im\alpha_k+i\alpha s)\rho}+O\left(s^2\right)\right) \\
        &=e^{-isz}\left(V(i\alpha s)+s\sum_{k=1}^{n}\sum_{m=1}^{\infty}\frac{\omega_k^me^{-imz_k}V(im\alpha_k+i\alpha s)}{m}+O\left(s^2\right)\right).
        \end{split}
    \end{align}
    Since the function $V(im\alpha_k+i\alpha s)$ is holomorphic at $s=0$, we have $V(im\alpha_k+i \alpha s)=V(im\alpha_k)+O(s)$ and
    \begin{align} \label{ExpandingV2}
        \sum_{m=1}^{\infty}\frac{\omega_k^me^{-imz_k}V(im\alpha_k+i\alpha s)}{m}=\sum_{m=1}^{\infty}\frac{\omega_k^me^{-imz_k}V(im\alpha_k)}{m}+O(s).
    \end{align}
    For brevity, defining the function
    \begin{align} \label{AuxFunctionDef2}
        \widetilde{f}(k)=\sum_{m=1}^{\infty}\frac{\omega_k^me^{-imz_k}V(im\alpha_k)}{m}
    \end{align}
    and substituting equation \eqref{ExpandingV2} in \eqref{L2SeriesToVfunction}, we obtain
    \begin{align} \label{L2ClosedFormFinally}
        \widetilde{L}(s;\underline{z},\underline{\alpha},\underline{\omega})=e^{-isz}\left(V(i\alpha s)+s\sum_{k=1}^{n}\widetilde{f}(k)+O\left(s^2\right)\right).
    \end{align}
    Taking linear terms in equation \eqref{L2ClosedFormFinally} results
    \begin{align*} 
        \begin{split}
            \underset{s=0}{\mathcal{LT}}\ \widetilde{L}(s;\underline{z},\underline{\alpha},\underline{\omega})&=\underset{s=0}{\mathcal{LT}}\left[e^{-isz}\left(V(i\alpha s)+s\sum_{k=1}^{n}\widetilde{f}(k)+O\left(s^2\right)\right)\right]\\
            &=\underset{s=0}{\mathcal{LT}}\left(e^{s(i\alpha/2-iz)}\phi(\alpha s)\right)+\sum_{k=1}^{n}\widetilde{f}(k).
        \end{split}
    \end{align*}
    By Lemma \ref{MeroPhiLemma}, we trivially have
    \begin{align*}
        \underset{s=0}{\mathcal{LT}}\left(e^{is(\alpha/2-iz)}\phi(\alpha s)\right)&=\underset{s=0}{\text{Res}}\frac{e^{s(i\alpha/2-iz)}\varphi_\alpha(s)}{s^2}\\
        &=\frac{c_{-1}}{2}\left(\frac{i\alpha}{2}-iz\right)^2+c_0\left(\frac{i\alpha}{2}-iz\right)+c_1=:\widetilde{F}(\underline{z};\underline{\alpha}).
    \end{align*}
    By the definition of the ddotted product, we get
    \begin{align} \label{LinearTermDDotted2}
    \begin{split}
        \widetilde{S}(\underline{z};\underline{\alpha},\underline{\omega})=\exp\left(-\underset{s=0}{\mathcal{LT}}\ \widetilde{L}(s;\underline{z},\underline{\alpha},\underline{\omega})\right)&=\exp\left(-\widetilde{F}(\underline{z};\underline{\alpha})-\sum_{k=1}^{n}\widetilde{f}(k)\right)\\
        &=e^{-\widetilde{F}(\underline{z};\underline{\alpha})}\prod_{k=1}^{n}\exp\left(-\widetilde{f}(k)\right).
    \end{split}
    \end{align}
    Exponentiating both sides of equation \eqref{AuxFunctionDef2} yields
    \begin{align*}
        \exp\left(-\widetilde{f}(k)\right)&=\exp\left(-\sum_{m=1}^{\infty}\frac{\omega_k^me^{-imz_k}V(im\alpha_k)}{m}\right)\\
        &=\exp\left(-\sum_{m=1}^{\infty}\sum_{\mathfrak{Im}(\rho)>0}\frac{-\omega_k^me^{im(\alpha_k\rho-z_k)}}{m}\right).
    \end{align*}
    We now interchange orders of summations as follows,
    \begin{align} \label{Expf(k)2}
        \begin{split}
            \exp\left(-\widetilde{f}(k)\right)&=\exp\left(-\sum_{\mathfrak{Im}(\rho)>0}\sum_{m=1}^{\infty}\frac{\omega_k^me^{im(\alpha_k\rho-z_k)}}{m}\right)\\
            &=\exp\left(\sum_{\mathfrak{Im}(\rho)>0} \log \left(1-\omega_ke^{i(\alpha_k\rho-z_k)}\right)\right)\\
            &=\prod_{\mathfrak{Im}(\rho)>0}1-\omega_ke^{i(\alpha_k\rho-z_k)}\\
            &=\left(\omega_ke^{-iz_k};e^{-i\alpha_k}\right)_\zeta.
        \end{split}
    \end{align}
    Finally, substituting equation \eqref{Expf(k)2} in \eqref{LinearTermDDotted2}, we have
    \begin{align*}
        \begin{split}
            \widetilde{S}(\underline{z};\underline{\alpha},\underline{\omega})&=e^{-\widetilde{F}(\underline{z};\underline{\alpha})}\prod_{k=1}^{n}\prod_{\mathfrak{Im}(\rho)>0}1-\omega_ke^{i(\alpha_k\rho-z_k)} \\
            &=e^{-\widetilde{F}(\underline{z};\underline{\alpha})}\prod_{k=1}^{n}\left(\omega_ke^{-iz_k};e^{-i\alpha_k}\right)_\zeta.
        \end{split}
    \end{align*}
    This completes the proof.
\end{proof}
\begin{corollary}
    For complex numbers $\alpha,z,\omega$ such that $2\pi>\mathfrak{Re}(z)\ge0\ge\mathfrak{Im}(z)$, $\alpha>0$ and $1\ge\left| \omega \right|$, the following equality holds.
    \begin{align*}
         \DDottedprod_{\mathfrak{Im}(\rho)>0}\left(e^{-i(\alpha\rho-z)}-\omega\right)=e^{-\widetilde{F}(z;\alpha)}\left(\omega e^{-iz};e^{-i\alpha}\right)_\zeta.
    \end{align*}
    Here $\widetilde{F}(z;\alpha)$ is a polynomial given by
     \begin{align*}
         \widetilde{F}(z;\alpha)=-\frac{\gamma+\log2\pi \alpha}{4\pi \alpha}\left(\frac{i\alpha}{2}-iz\right)^2+\frac{7(i\alpha-2iz)}{16}+c_1.
     \end{align*}
\end{corollary}
Similarly, we can calculate the discrepancy of the sequences $e^{-i(\alpha_k\rho-z_k)}-\omega_k$ as follows with the help of Theorem \ref{Mainthm2}.
\begin{corollary}
    The discrepancy $\widetilde{F}$ of the sequences $e^{-i(\alpha_k\rho-z_k)}-\omega_k$, $k=1,\ldots,n$ is given by
    \begin{align*}
        \widetilde{F}=\sum_{k=1}^{m}\widetilde{F}(z_k;\alpha_k)-\widetilde{F}(\underline{z};\underline{\alpha}).
    \end{align*}
\end{corollary}
The discrepancy between sequences of form $\prod_{j=1}^{m_k} \left(e^{-i(\alpha_{j,k}\rho-z_{j,k})}-\omega_{j,k}\right)$ may also be calculated in an analogous manner.

\section{Concluding Remarks}

We conclude our paper with an indication of what could be pursued in further research. The following relationship with the Riemann hypothesis and zeta regularized products is established in \cite{WakayamaRemarks}.
\begin{theorem}\label{RHThm}
    Suppose that the regularized product $\DDottedprod_{\mathfrak{Im}(\rho)>0}|\sin\alpha(\rho-x)|$ exists. Then, the following two conditions are equivalent.
    \begin{enumerate}
        \item The equality 
            \begin{align*}
                \DDottedprod_{\mathfrak{Im}(\rho)>0}|\sin\alpha(\rho-x)|=e^{-R_\alpha(\mathfrak{Re}(x))}\left|\DDottedprod_{\mathfrak{Im}(\rho)>0}\sin\alpha(\rho-x) \right|
            \end{align*}
            holds for two distinct values of $\alpha>0$ where $\mathfrak{Im}(x)\le0$.
            \item An asymptotic formula $\sum_{\mathfrak{Im}(\rho)>0}(\mathfrak{Re}(\rho)-1/2)^2 e^{-t\mathfrak{Im}(\rho)}=O(\log t)$ holds as $t\to 0$.
    \end{enumerate}
    Here $R_\alpha(x)$ is a quadratic polynomial defined by
    \begin{align*}
        R_\alpha(x)=-\frac{\gamma+\log2\pi \alpha}{4\pi \alpha}\left( \alpha\left( x-\frac{1}{2} \right)+\frac{\pi}{2} \right)^2.
    \end{align*}
\end{theorem}
It is furthermore remarked that if the series
\begin{align*}
    \sum_{\mathfrak{Re}(\tau)>0}(\mathfrak{Re}(\tau))^x (\mathfrak{Im}(\tau))^2
\end{align*}
converges for all $x>0$, then the regularized product $\DDottedprod_{\mathfrak{Im}(\rho)>0}|\sin\alpha(\rho-x)|$ exists. The question of whether an analogue of the infinite product in Theorem \ref{firstthm} arises for $\DDottedprod_{\mathfrak{Im}(\rho)>0}|\sin\alpha(\rho-x)|$. It seems a natural direction to investigate whether theorems of this type also holds for generalized products presented in this paper.\\
\newline

Another direction of research is to examine such products for other zeta functions. A similar investigation of such trigonometric products have been conducted for the case of Selberg zeta functions in \cite{WakayamaRemarks}. Furthermore, a similar criterion of above type is shown for Selberg's $1/4$ conjecture for congruence subgroups of $PSL_2(\mathbb{R})$. The criterion is much stronger, in the sense that it only needs such an equation to hold for a single value of $\alpha>0$ and is directly equivalent to the Selberg's $1/4$ conjecture. The authors are currently working on generalizing the results of this paper to suitable Selberg zeta functions. It would be interesting to study how the generalized products relate to Selberg's $1/4$ conjecture in the future.

\bibliographystyle{abbrvnat}

\end{document}